\documentclass[11pt,reqno]{amsart}
\usepackage{amssymb,amsmath,amsthm,amsfonts}
\usepackage{graphicx}
\usepackage{xcolor}

\usepackage[capitalize,noabbrev,nameinlink]{cleveref}

\crefformat{equation}{(#2#1#3)}
\crefformat{figure}{#2Figure #1#3}
\crefformat{enumi}{(#2#1#3)}
\crefrangeformat{enumi}{(#3#1#4)--(#5#2#6)}
\crefformat{enumii}{(#2#1#3)}
\crefformat{enumiii}{(#2#1#3)}
\crefformat{theorem}{#2Theorem~#1#3}
\crefformat{lemma}{#2Lemma~#1#3}
\crefformat{proposition}{#2Proposition~#1#3}
\crefformat{corollary}{#2Corollary~#1#3}
\crefformat{claim}{#2Claim~#1#3}
\crefformat{section}{#2Section~#1#3}
\crefformat{subsection}{#2Section~#1#3}

\usepackage[numbers,sort]{natbib}

\DeclareMathOperator{\diam}{diam}

\DeclareMathOperator{\adim}{dim_A}
\newcommand{\N}{\mathbb{N}}                      
\newcommand{\R}{\mathbb{R}}                      

\newcommand{\aspecdim}[1][]{\operatorname{\overline{dim}_A^{\, #1}}}

\newcommand{\de}{\delta}							
\newcommand{\la}{\lambda}
\newcommand{\eps}{\epsilon}							

\newcommand{\ovdimB}{{\overline{\dim_B}}}

\theoremstyle{plain}
\newtheorem{theorem}{Theorem}
\newtheorem{theoremA}{Theorem}

\newtheorem{lemma}[theorem]{Lemma}

\theoremstyle{definition}

\numberwithin{theorem}{section} \numberwithin{equation}{section}

\title{Minkowski weak embedding theorem}

\author{Efstathios-K. Chrontsios-Garitsis}
\address[Efstathios-K. Chrontsios-Garitsis]{Department of Mathematics, University of Tennessee,
	Knoxville, 1403 Circle Dr, Knoxville, TN 37966, USA}
\email{echronts@utk.edu}
\urladdr{https://sites.google.com/view/efstathioschrontsiosgaritsis}

\author{Sascha Troscheit}
\address[Sascha Troscheit]{Mathematical Sciences Research Unit, PO Box 8000, FI-90014 University of Oulu, Finland}
\email{arxiv@troscheit.eu}
\urladdr{https://www.troscheit.eu}

\begin{document}
	\maketitle

	\begin{abstract}
		A well-known theorem of Assouad states that metric spaces satisfying the doubling property can be
		snowflaked and bi-Lipschitz embedded into Euclidean spaces. Due to the invariance of many geometric
		properties under bi-Lipschitz maps, this result greatly facilitates the study of such spaces. We
		prove a non-injective analog of this embedding theorem for spaces of finite Minkowski dimension.
		This allows for non-doubling spaces to be weakly embedded and studied in the usual Euclidean
		setting. Such spaces often arise in the context of random geometry and mathematical physics with
		the Brownian continuum tree and Liouville quantum gravity metrics being prominent examples. 
	\end{abstract}
	
	\section{Introduction}
	Doubling spaces constitute one of the most studied classes of metric spaces within many areas of
	mathematics. This is to a great extent due to the Assouad embedding theorem \cite{Assouad83}, a
	celebrated result with applications in analysis, geometry and geometric group theory. See \cite{Assouad_Analysis}, \cite{Assouad_geom}, \cite{Assouad_GeomGroupTh} for a non-exhaustive list of applications in the aforementioned areas, as well as
	\cite{hei:lectures} and the references therein. Assouad's theorem loosely states that any doubling metric space can be embedded into a Euclidean space $\R^M$ for some $M \in \N$
	using a bi-Lipschitz map, up to some arbitrarily small distortion of the space. The result was later quantitatively improved by Naor
	and Neiman \cite{Naor_Neiman_Assouad}, while more recently David and Snipes \cite{Assouad_David_Snipes} gave a simplified and
	constructive proof of the following stronger version of the theorem.
	\begin{theoremA}
		Suppose $(X,d)$ is a doubling metric space with doubling constant $C\geq1$.
		There exist $M\in\N$ and, for all $\tfrac12 < \eps < 1$, a constant $K=K(C,\eps)$ and a mapping
		$F:X\to\R^{M}$ such that 
		\[
		K^{-1} d(x,y)^{\eps} \leq |F(x) - F(y)| \leq K d(x,y)^{\eps}
		\]
		for all $x,y\in X$.
	\end{theoremA}
	Note that the exponent $\eps$ quantifies the amount of distortion we have to impose on $X$ before we
	bi-Lipschitz embed the space in a Euclidean space. Such distortion is called
	\textit{$\eps$-snowflaking} of the metric space $(X,d)$, and the resulting metric space $(X,d^\eps)$
	is often denoted by $X^\eps$. While $\eps$ can be chosen to be as close to $1$ as needed, there are
	cases, such as the Heisenberg group, where a bi-Lipschitz embedding is not possible for $\eps=1$,
	see  \cite[p.~99]{hei:lectures}.
	
	Related to the Assouad embedding result is the notion of the \textit{Assouad dimension}, 
	which can be interpreted as a quantisation of the doubling property for metric spaces.
	Given a metric space $(X,d)$, the Assouad dimension of $X$ is defined by
	\begin{equation} \label{eq:adimdef}
		\adim E = \inf\left\{ s>0 \mid (\exists K>0)(\forall 0<r<R) \,\,\, \sup_{x\in X}N_r(B(x,R)) \leq K\left(
		\frac{R}{r} \right)^{s}\right\},
	\end{equation}
	where $N_r(B(x,R))$ is the minimal number of balls of radius $r$ required to cover $B(x,R)$. The
	Assouad dimension has received a lot of attention in recent years due to its many applications in
	fractal geometry, dynamics, harmonic analysis, number theory, and other fields; see 
	\cite{wang2024assouaddimensionkakeyasets,BaranyBook,Anderson21,Saito20} for some recent examples as well as the text book \cite{Fraser2020} for a comprehensive exposition.
	There is a close relation of the Assouad dimension notion with the doubling property. Namely, a metric space $(X,d)$ is doubling if, and only if,
	$\adim X <\infty$. This property enables us to state Assouad's theorem for spaces with finite
	Assouad dimension and motivates the investigation of similar embedding results for other dimension
	notions, which are smaller than the Assouad dimension. A recent example is a non-injective embedding result by David \cite{Guy_Nagata}, which
	considers the Nagata dimension.
	
	The Minkowski dimension, or box-counting dimension, of a metric space $X$ is a classical notion of
	dimension that has been popular
	within many fields of research.
	Some examples include its use in studying partial differential equations  \cite{MinkNavierStokesEq}, signal processing \cite{MinkSignal}, and mathematical physics \cite{MinkPhysics}.
	Its intuitive definition through ``intersection with grid squares'' (see \cite[p.~6, 8]{Fraser2020}) makes it easy to approximate
	using numerical methods and many articles in applied fields simply refer to the Minkowski dimension
	as \textit{the} ``fractal dimension".
	For a bounded metric space $(X,d)$ the (upper) Minkowski dimension $\ovdimB X$
	can be defined in a similar way to \cref{eq:adimdef}:
	\begin{equation}\label{eq:boxdim}
		\ovdimB X = \inf\left\{ s>0 \mid (\exists K>0)(\forall 0<r<\diam X) \,\,\, \sup_{x\in X} N_r(X) \leq K\left(
		\frac{1}{r} \right)^{s}\right\}.
	\end{equation}
	We make no reference to the lower Minkowski dimension and, hence, we omit the adjective ``upper" in what follows. From the definitions \cref{eq:adimdef}, \cref{eq:boxdim} it is immediate that $\ovdimB X \leq \adim X$, by fixing $R = \diam (X)$. 
	
	In analogy to finite Assouad dimension implying doubling and, hence, the Assouad embedding theorem, one
	may ask whether a similar embeddability result is guaranteed using the weaker assumption that $\ovdimB X
	<\infty$. In this article we prove such a qualitative embedding result.
	\begin{theorem}\label{thm:main}
		Let $(X,d)$ be a compact metric space with $\ovdimB X<\infty$, and $\eps\in(1/2,1)$. There are $\tau\in (0,1)$ and $n_0\in \N$ such that for every $n> n_0$  there is an $L$-Lipschitz map $F:X^\eps\rightarrow\R^{M}$, where $L=L_{\eps,n}$ and $M=M_{n}$, such that $L^{-1}d(x,y)^\eps\leq |F(x)-F(y)|$ for all $x$, $y\in X$ for which $d(x,y)\geq 4\tau^{2n}$. 
	\end{theorem}
	We remark that the lower Lipschitz constant of the map in \cref{thm:main} can in fact be taken to be independent of $n$ (see \eqref{eq: F lower-Lip}).
	
	It should be noted that there is a large class of metric
	spaces arising in probability which are not doubling, but have finite Minkowski dimension. For
	instance, the Brownian continuum random tree has
	Hausdorff and Minkowski dimension equal to $2$, see
	\cite{Continuum_random_tree_I,Continuum_random_tree_II,Continuum_random_tree_III}; whereas the Brownian map has
	Hausdorff and Minkowski dimension $4$, see \cite{Brownian_map}.
	However, the Assouad dimension in these examples is infinite and bi-Lipschitz, or even
	quasi-symmetric, embeddings into Euclidean spaces do not exist, even after suitable snowflaking, see
	\cite{Sascha_non-doubling_Minkowski_finite}. The same is true for more general random metric spaces
	appearing in mathematical physics, such as in Liouville quantum gravity metrics, see \cite{Liouville_quantum_gravity}. 
	
	\medskip
	\textbf{Organisation.} This paper is organised as follows. In \cref{sec:background} we establish the notation and outline the
	necessary background and tools from fractal geometry. In \cref{sec: proof_main} we construct a
	mapping that satisfies the properties needed to prove \cref{thm:main}.
	\cref{sec:questions} contains some remarks and potential future directions.
	
	\medskip
	\textbf{Acknowledgments.} 
	This manuscript was finalised while ST was visiting the Fields Institute as part of the
	\textit{Randomness and Geometry} Semester Programme. ST is grateful for the hospitality and
	productive atmosphere provided at the Fields Institute. 
	ST's research was supported by the European Research Council Marie Curie-Sk\l{}odowska Personal Fellowship \#101064701.
	
	\section{Background}\label{sec:background}
	
	Given a metric space $(X,d)$, for $x\in X$ and $r>0$ we denote by $B(x, r)$ the closed ball with center $x$ and radius $r$. For $M\in \N$ and $L\geq 1$, we say the map $F:X\rightarrow \R^M$ is \textit{$L$-Lipschitz} if
	$$
	|F(x)-F(y)|\leq L \, d(x,y),
	$$ for all $x, y\in X$. Moreover, if $L^{-1} d(x,y)\leq |F(x)-F(y)|$ also holds for all $x,y \in X$,
	we say that $F$ is \textit{($L$-)bi-Lipschitz}. We denote the \textit{Lipschitz norm} of $F$ by
	$$
	||F||_{\text{Lip}} := \sup\limits_{x, y\in X, \,\, x\neq y}\frac{|F(x)-F(y)|}{d(x,y)}.
	$$

	Our proof relies on the relation between the Minkowski dimension and the Assouad spectrum. The
	latter is a continuum of dimensions which aim to interpolate between the Minkowski dimension and the
	Assouad dimension. The notion was first introduced by Fraser and Yu in \cite{fy:assouad-spectrum},
	but has received subsequent interest as a means to determine exact relative scaling behaviour, see
	e.g.~\cite{Banaji24,Banaji23,Fraser24,Roos23,ChronTysQCspec,Chron24,Hare22}.
	Given $\theta\in (0,1)$, the \textit{(regularised) $\theta$-Assouad spectrum} of $X$ is defined, in analogy to \cref{eq:adimdef} and \cref{eq:boxdim}, as
	\begin{equation}\label{eq: def of assouad spec}
		\aspecdim[\theta] X = \inf\left\{ s>0 \mid (\exists K>0)(\forall 0<r\leq R^{1/\theta}<R<1)
		\;\sup_{x\in X}N_r(B(x,r^{\theta}) \leq K\left(
		\frac{R}{r} \right)^{s}\right\}.
	\end{equation}
	The term ``upper" Assouad spectrum has also been used in the literature, see
	\cite[p.~5]{ChronTysQCspec}, \cite[Section 3.3]{Fraser2020}, and \cite{Fraser19}. For a fixed space $X$, the spectrum $\aspecdim[\theta]X:(0,1) \to [0,\infty]$ is continuous and increasing in
	$\theta$. It continuously tends to the Minkowski dimension as $\theta\rightarrow 0^+$ and typically
	(but not necessarily) approaches the Assouad dimension as $\theta\to1^-$, see
	\cite{Fraser2020,Rutar24,Banaji23} for a detailed treatment of its behaviour.

	Recall that a metric space $(X,d)$ is called \textit{doubling} with doubling constant $C_0>0$, if
	for any $R>0$, every ball of radius $R$ can be covered by at most $C_0$ balls of radius $R/2$. Note
	that for fixed $\lambda > 1$, every ball of radius $\lambda R$ in $X$ can be covered by at most
	$C_0^{\lambda \log_2 C_0}$ balls of radius $R$.
	For $M\in \N$, we consider $\R^M$ equipped with its
	Euclidean metric and denote the doubling constant of $\R^M$ by $C_M\geq 1$.
	
	In \cref{sec: proof_main} we take advantage of the fact that $\ovdimB X <\infty$ is equivalent to the existence of some
	$\theta\in(0,1)$ such that $\aspecdim[\theta]X <\infty$, see \cite[Theorems 3.3.1 and 3.3.6]{Fraser2020}.
	This allows us to use the finiteness of the Assouad spectrum to achieve a weaker doubling-like condition for our space, and adjust the proof of David and Snipes from \cite{Assouad_David_Snipes} accordingly.
	As mentioned in the introduction, the Assouad dimension of $X$ is finite if, and only if, $X$ is doubling. On the other hand,  if
	$\aspecdim[\theta] X<\de<\infty$ for some $\theta\in (0,1)$, then, by choosing $r=R^{1/\theta}$ in \cref{eq: def of assouad spec}, every ball of radius $\la R\in (0,1)$
	can be covered by at most $\la^\delta C_{\theta,\delta} R^{(1-1/\theta)\de}$ balls of radius $R^{1/\theta}$.
	Similarly to Assouad's proof, this is a core property of $X$ for the proof of \cref{thm:main}. We
	say that a metric space $X$ satisfying this property is \textit{$(\theta, \delta)$-quasidoubling}
	with quasidoubling constant $C_{\theta,\de}$, cf.~the notion of quasidoubling for measures in
	\cite{Hare20}.

	\section{Constructing the embedding}\label{sec: proof_main}
	Let $(X,d)$ be a compact metric space with $\ovdimB X<\infty$. Then $d_\theta:=\dim_A^\theta X<\infty$ for some $\theta\in (0,1)$. In fact, a closer analysis of Theorem 3.3.1 in \cite{Fraser2020} reveals that the variable $\theta$ can be chosen arbitrarily close to $1$. Fix $\de>\dim_A^\theta X$, $\eps\in (1/2,1)$ and set the quasidoubling constant of $X$ to be $C=C_{\theta,\de}$ to ease notation. We can scale the metric of $X$ in a bi-Lipschitz way, which does not change the Minkowski dimension of $X$ (see Table 2.1 in \cite{Fraser2020}). Hence, without loss of generality we may assume that $\diam X<1$.

	The proof of \cref{thm:main} requires some preparation. Fix $\tau \in(0,1/2)$ small enough for all the following inequalities to hold:
	\begin{equation}\label{eq: tau0}
		\tau<\min\{ 1/2^{2\eps},1-\eps \}
	\end{equation}
	\begin{equation}\label{eq: tau1}
		\tau^{2(\theta-1)}>2
	\end{equation}
	\begin{equation}\label{eq: tau2}
		e^{-2\tau \log(1/\tau)}\leq 1-\tau\log (1/\tau)
	\end{equation}
	\begin{equation}\label{eq: tau3}
		1-4 \tau^{2\eps-1}\geq 1/2
	\end{equation}
	\begin{equation}\label{eq: tau4}
		-7\log (7\tau)> \max \{ \log (40^\de C^2),\, \log 21 \}.
	\end{equation}
	While it is certainly not clear at this point why these specific conditions are chosen for $\tau$, their need is justified in calculations that follow.
	
	Set 
	\begin{equation}\label{eq:rk_scales_def}
		r_k := \tau^{2k}, \,\,\, k \in \mathbb{Z}
	\end{equation} to be the scales we use to construct our map, and $n_0$ to be the largest integer such that $r_{n_0} \geq \text{diam}(X)$.
	
	Fix some integer $n>n_0$ for what follows. For each integer $k \in[n_0,n]$, fix a maximal collection $X_k:=\{x_j^k\}_{j \in J_k}$ of points in $X$ with $d(x_i^k, x_j^k) \geq r_k$ for all distinct $i, j \in J_k$. By maximality of this collection we have
	\begin{equation}\label{eq:X_covering_def}
		X \subset \bigcup_{j \in J_k} B(x_j^k, r_k).
	\end{equation}
	
	We write $B_j^k=B(x_j^k,r_k)$, and $x_j=x_j^k$ when the level $k$ is clear. By the assumption that $X$ is compact, we have that $J_k$ is finite. 
	
	For $x\in X$, set $N_k(x)$ to be the number of points $x_j$, $j \in J_k$, for which $d(x_j, x) \leq r_k^\theta$. Recall that by the $(\theta,\delta)$-quasidoubling property of $X$, we can cover $B(x, r_k^\theta)$ with at most $C \left( \frac{r_k^\theta}{3^{-1}r_k}\right)^\de$ balls of radius $r_k/3$, which we denote by $B_{x,\ell}$, $\ell\in J_k^x$. Due to $d(x_i, x_j) \geq r_k$ for $i \neq j$, every $B_{x,\ell}$ contains at most one point $x_j$. But all the points $x_j$ that lie in $B(x, r_k^\theta)$ are contained in some $B_{x,\ell}$. Therefore, we have
	\begin{equation}\label{eq:N(x)}
		N_k(x) \leq C \left( \frac{r_k^\theta}{3^{-1}r_k}\right)^\de=3^\delta C r_k^{\de(\theta-1)}
	\end{equation}
	for all $x \in X$.
	
	Following the idea of Assouad's proof in \cite{Assouad83}, but now with the bound \eqref{eq:N(x)}
	quantitatively depending on the level $k$, set $\Xi_n = \{1, 2, \ldots, [3^\de C
	r_n^{\de(\theta-1)}]\}$ to be a set of colors at level $n\geq k$, where $N_n:=[3^\delta C
	r_n^{\delta(\theta-1)}]\in \N$  and $[.]$ denotes the ceiling function. Since $n$ is fixed, we write $\Xi=\Xi_n$ for what follows.  Fix an enumeration of $J_k$, and define $\xi:J_k\rightarrow\Xi$ inductively with $\xi(1)=1$ and
	$$
	\xi(j)= \min\left\{ \xi\in\Xi: \xi\neq \xi(j') \,\, \text{for all}\,\, j'\in \{1, \dots, j\} \,\, \text{with} \,\, x_{j'}\in B(x_j,r_k^\theta)  \right\},
	$$for all $j\in J_k$. Namely, $\xi(j)$ is the smallest color not already assigned to an earlier point of $X_k$ that is close to $x_j$. Note that the function $\xi$ is well defined due to \eqref{eq:N(x)} and $k\leq n$. We also have by definition that
	$\xi(i) \neq \xi(j) \text{ for all distinct } i, j \in J_k \text{ such that }  d(x_i, x_j) \leq r_k^\theta$.
	Moreover, for all $\xi \in \Xi$ set $J_k(\xi) := \{j \in J_k \mid \xi(j) = \xi\}$, which implies
	\begin{equation}\label{eq:same_color_points_far}
		d(x_i, x_j) > r_k^\theta \text{ for } i, j \in J_k(\xi) \text{ such that } i \neq j.
	\end{equation}
	
	For each $j \in J_k$, 
	define $\phi_j:=\phi_j^k: X \rightarrow \R$ by $\phi_j (x) = \max\{0, 1 - r_k^{-1}\operatorname{dist}(x, B_j)\}$ for all $x\in X$. This choice ensures that
	$0 \leq \phi_j (x) \leq 1$ for all $x\in X$ and $\phi_j (x) = 1$ for all $x\in B_j^k$.
	Moreover, due to this choice, the following conditions hold:
	\begin{equation}\label{eq:phi_0_outside_2B}
		\phi_j (x) = 0 \text{ for } x \in X \setminus 2B_j^k,
	\end{equation}
	and
	\begin{equation}\label{eq: phi_j is Lip}
		\|\phi_j\|_{\text{Lip}} \leq r_k^{-1}.
	\end{equation}
	
	For each $\xi \in \Xi$, we construct a coordinate mapping $F^\xi : X \to \mathbb{R}^M$ and a
	modified version $\widetilde{F}^{\xi} : X \to \mathbb{R}^M$, where $M=M_{n}$  is a large enough integer
	we pick appropriately depending only on $\theta$, $\de$ and $n$ (see \eqref{eq: M_n choice}). We
	then combine these two mappings by taking their direct product for all colors in $\Xi$ to get $F : X
	\to \R^{\widetilde{M}}$, for $\widetilde{M}=2M_n N_n$ (recall that $N_n$ is the number of colors).
	
	We first make preparations for the construction of the coordinate maps $F^\xi$ inductively. Set
	\begin{equation}\label{eq: F ksi def}
		F^\xi(x) = \sum_{n_0\leq k \leq n} r_k^\eps f^{\xi}_k(x):X\rightarrow\R^{M},
	\end{equation}
	where
	\begin{equation}\label{eq: f^ksi def}
		f^{\xi}_k(x) = \sum_{j \in J_k(\xi)} v_j^k \phi_j(x):X\rightarrow\R^{M}
	\end{equation}
	for all $k\in[n_0,n]$, and with vectors $v_j^k \in \mathbb{R}^{M_k}$ chosen inductively later, in
	\cref{lem: choice of v_j}, so that
	\begin{equation}\label{eq: v_j in ball}
		v_j^k \in B(0, \tau^2) \subset \mathbb{R}^{M}.
	\end{equation}
	We now outline certain necessary properties of such maps, which hold irrespectively of the specific choice of the vectors $v_j^k$, other than the bound on their magnitude by $\tau^2$. By \eqref{eq: tau1} and \eqref{eq: v_j in ball} we get
	$$
	||f_\xi^k||_\infty \leq \tau^2,
	$$
	since all maps $\phi_j$, $j \in J_k(\xi)$, have disjoint supports due to $d(x_i,x_j)>r_k^\theta>2r_k$, \eqref{eq:same_color_points_far} and \eqref{eq:phi_0_outside_2B}. Hence, the series $F^\xi$ in \eqref{eq: F ksi def} converges. Moreover, for $k<n$ we set
	\begin{equation}\label{eq: F_k^ksi def}
		F_k^\xi(x) = \sum_{n_0 \leq \ell \leq k} r_\ell^\eps f_\ell^\xi(x),
	\end{equation}
	which implies that
	\begin{equation}\label{eq: F_k^ksi convergence}
		||F^\xi - F^\xi_k||_\infty \leq \sum_{\ell > k} r^\eps_\ell \tau^2 = r^\eps_{k+1} \tau^2 \sum_{\ell \geq 0} \tau^{2\ell\eps} \leq 2\tau^2 r^\eps_{k+1},
	\end{equation}
	due to $r_\ell = \tau^{2\ell}$ and by \eqref{eq: tau0}. We also have 
	$||f^\xi_k||_{\text{Lip}} \leq \tau^2 r_k^{-1}$
	by \eqref{eq: phi_j is Lip} and due to $\phi_j=\phi_j^k$ being supported in disjoint balls $2B_j^k$. Summing over all levels up to $k$ and using \eqref{eq: F_k^ksi def} results in the following:
	\begin{equation}\label{eq: F_k^xi Lip norm pre}
		||F^\xi_k||_{\text{Lip}} \leq \sum_{\ell \leq k} r^\eps_\ell ||f^\xi_\ell||_{\text{Lip}} \leq \tau^2 \sum_{\ell \leq k} r^{\eps-1}_{\ell} = \tau^2 r^{\eps-1}_{k} \sum_{\ell \leq k} \tau^{2(\ell-k)(\eps-1)} = \tau^2 r_k^{\eps-1} (1 - \tau^{2(1-\eps)})^{-1}.
	\end{equation}
	Since $\tau \leq 1 - \eps$ by \eqref{eq: tau0}, we have
	$$
	\ln(\tau^{2(1-\eps)}) = 2(1 - \eps) \ln(\tau) = -2(1 - \eps) \ln\left(\frac{1}{\tau}\right) \leq -2\tau \ln\left(\frac{1}{\tau}\right),
	$$
	which implies $\tau^{2(\eps-1)} \leq e^{-2\tau \ln(1/\tau)} \leq (1-\tau)\ln(1/\tau)$ by \eqref{eq: tau2}. Thus, using this on \eqref{eq: F_k^xi Lip norm pre} we get
	\begin{equation}\label{eq: F_k^ksi Lip norm}
		||F^\xi_k||_{\text{Lip}} \leq \frac{\tau}{\ln(1/\tau)} r_k^{\eps-1} \leq r_k^{\eps-1}.
	\end{equation}
	
	We now focus on the inductive construction of the vectors $v_j^k$, $j \in J_k$. Fix an integer $k \in [n_0,n]$, a color $\xi \in \Xi$, and suppose that the map $F^\xi_{k-1}$ is constructed, with the trivial base case $F_{n_0-1}=0$. Fix an increasing order on $J_k(\xi)$ with respect to which we construct $F^\xi_k$. 
	Recall that
	\begin{equation}\label{eq: Fk is F_k-1 + phis}	
		F^\xi_k(y) = F^\xi_{k-1}(y) + r^\eps_k f^\xi_k(y) = F^\xi_{k-1}(y) + r^\eps_k \sum_{i\in J_k(\xi)} v_i^k \phi_i(y),
	\end{equation}
	and for each $j \in J_k(\xi)$ define the partial sum $G^\xi_{k,j}$ such that
	\begin{equation}\label{eq: G_k def}
		G^\xi_{k,j}(y) = F^\xi_{k-1}(y) + r^\eps_k \sum_{i\in J_k(\xi) ; i<j} v_i^k \phi_i(y),
	\end{equation} for all $y\in X$.
	
	\begin{lemma}\label{lem: choice of v_j}
		There is some integer $M=M_{n}$ such that for each $j \in J_k(\xi)$, there are vectors $v_j^k \in B(0, \tau^2)\subset \R^M$ so that
		\begin{equation}\label{eq: Lower Lip Lemma}
			|F^\xi_k(x) - G^\xi_{k,j}(y)| \geq \tau^3 r^\eps_k,
		\end{equation}
		for all $x \in B_j^k$, $y \in B(x_j, 10\tau^{-2} r_k) \setminus 2B_j^k$.
	\end{lemma}
	\begin{proof} Fix $j\in J_k(\xi)$. Since for what follows $k$ is also fixed, set $B_j:=B_j^k$ and $v_j := v_j^k$ to simplify notation. Note that $\phi_i$ is supported in $2B_i$ by \eqref{eq:phi_0_outside_2B}, and $2B_i$ never meets $B_j$ for $i\neq j$ by \eqref{eq:same_color_points_far}. These imply that for all $x \in B_j$, we have $\phi_j(x) = 1$ and $\phi_i(x)=0$ for all $i\neq j$. Thus, using \eqref{eq: F ksi def} and \eqref{eq: f^ksi def}, we get
		\begin{equation}\label{eq: F_k is F_k-1+..}
			F^\xi_k(x) = F^\xi_{k-1}(x) + r^\eps_k f^\xi_k(x) = F^\xi_{k-1}(x) + r^\eps_k v_j.
		\end{equation}
		Following the proof of \eqref{eq: F_k^ksi Lip norm}, by adding fewer terms, we also have
		$$
		||G^\xi_{k,j}||_{\text{Lip}} \leq r_k^{\eps-1}.
		$$
		
		The above allows us to replace $B_j$ and $B(x_j, 10\tau^{-2} r_k) \setminus 2B_j$ with discrete sets. Set $\eta = \tau^3 r_k^{1/\theta}$, and pick an $\eta$-dense set $E_1$ in $B_j$, and an $\eta$-dense set $E_2$ in $B(x_j, 10\tau^{-2} r_k) \setminus 2B_j$. It is enough to show that $v_j$ can be picked so that
		\begin{equation}\label{eq: F_k-G_k}
			|F^\xi_k(x') - G^\xi_{k,j}(y')| \geq 3\tau^3 r^\eps_k,
		\end{equation} for all $x'\in E_1$, $y'\in E_2$.
		
		Indeed, if $x \in B_j$, there is $x' \in E_1$ such that
		$$
		|F^\xi_k(x') - F^\xi_k(x)| \leq  ||F^\xi_k||_{\text{Lip}}\, d(x',x)  \leq r_k^{\eps-1} \tau^3 r_k^{1/\theta} \leq \tau^3 r^\eps_k,
		$$
		and, similarly, if $y \in B(x_j, 10\tau^{-2} r_k) \setminus 2B_j$, there is $y' \in E_2$ such that
		$$
		|G^\xi_{k,j}(y) - G^\xi_{k,j}(y')| \leq ||G^\xi_{k,j}||_{\text{Lip}}\, d(y,y') \leq \tau^3 r^\eps_k.
		$$
		The above two inequalities, along with \eqref{eq: Lower Lip Lemma}, imply \eqref{eq: F_k-G_k} and
		conclude the proof. 
		
		In order to construct the appropriate  $\eta$-dense sets $E_1\subset B_j$, $E_2\subset B(x_j, 10\tau^{-2} r_k) \setminus 2B_j$ for \eqref{eq: F_k-G_k} to hold, we also need to bound $|E_1|$ and $|E_2|$ accordingly. By $\de\in (\dim_A^\theta X, \infty)$, we can cover $B_j = B(x_j, r_k)$ by at most $2^\de C r_k^{(1-1/\theta)\de}\tau^{-3\de}$ balls $B_{j,\ell}$ of radius $\eta/2$, since $\eta/2=\tau^3 {r_k}^{1/\theta}/2<{r_k}^{1/\theta}$. Fix an element of $B_j$ in each such ball $B_{j,\ell}$ that intersects $B_j$. This collection of elements forms $E_1$, which is an $\eta$-dense subset of $B_j$, such that $$
		|E_1| \leq 2^\de C r_k^{(1-1/\theta)\de}\tau^{-3\de}.
		$$ 
		Similarly, we can find $E_2$, an $\eta$-dense subset of $B(x_j, 10\tau^{-2} r_k) \setminus 2B_j$, so that
		$$
		|E_2| \leq C \left(\frac{10\tau^{-2}r_k}{\tau^{3}r_k^{1/\theta}2^{-1}}\right)^{\de} = 20^\de C \tau^{-5\de}r_k^{(1-1/\theta)\de},
		$$
		since $\tau^3 {r_k}^{1/\theta}/2<{r_k}^{1/\theta}<(10\tau^{-2}r_k)^{1/\theta}$.
		Therefore, we need to pick suitable vectors $v_j$ and check \eqref{eq: F_k-G_k} for a total number 
		$$
		|E_1||E_2| \leq 40^\de C^2 \tau^{-8\de}r_k^{2(1-1/\theta)\de}
		$$ of pairs of points $(x', y')$ in $E_1\times E_2$. 
		
		We now choose the vectors $v_j$ from a maximal finite set $V$ lying in $B(0, \tau^2) \subset
		\mathbb{R}^{M_k}$, whose points are pair-wise at least $7\tau^3$ away from each other. Recall that
		$G^\xi_{k,j}$ does not depend on $v_j=v_j^k$ due to \eqref{eq: G_k def}. Therefore, for each pair
		$(x', y')\in E_1\times E_2$, different choices of $v_j \in V$ yield the same value of
		$G^\xi_{k,j}(y')$, and values of $F^\xi_k(x')$ that lie at a distance at least $7\tau^3 r^\eps_k$
		apart. To see the latter, suppose $x'$ and $\widetilde{x}'$ correspond to the choices $v_j$ and $u_j$ in $V$, respectively. Then by \eqref{eq: F_k is F_k-1+..} we have
		$$
		|F_k^\xi(x')-F_k^\xi(\widetilde{x}')|=|r_k^\eps(v_j-u_j)|\geq r_k^\eps 7 \tau^3,
		$$ since $v_j$ and $u_j$ lie at a distance at least $7\tau^3$ away from each other by definition of $V$.
		Hence, the inequality \eqref{eq: F_k-G_k} for a specific pair $(x', y')$ might fail for at most one choice of $v_j \in V$, which implies that it is enough to show that $|V|>|E_1||E_2|$. Recall that $C_M$ denotes the doubling constant of $\mathbb{R}^M$, where $M=M_n$ is a large enough integer we have not yet determined. By the doubling property of $\R^M$ and the fact that elements of $V$ lie at a pair-wise distance of at least $7\tau^3$ apart, we have $|V| \geq C_M (1/7\tau)^M$. For the right-hand side to be larger than $|E_1||E_2|$ it is enough to have 
		\begin{equation}\label{eq: M_n size}
			M>\frac{\log (40^\delta C^2 \tau^{-8\delta})+2(1-1/\theta)\delta \log r_k-\log C_M}{-7\log 7\tau}.
		\end{equation} Note that choosing $M$ to be equal to
		$$
		\left[\frac{\log (40^\delta C^2 \tau^{-8\delta})+2(1-1/\theta)\delta \log r_n}{-7\log 7\tau}\right]= \left[ \frac{\log 40^\de C^2-\log 7^{-8\de+4n\de (1-1/\theta)}}{-7 \log7\tau}+ \frac{-8\de+4n\de (1-1/\theta)}{-7} \right]
		$$ would be enough for \eqref{eq: M_n size} to hold true, but then $M$ would implicitly depend on $\eps$ through its dependence on $\tau$. For this reason, we choose
		\begin{equation}\label{eq: M_n choice}
			M:= \left[ \frac{6 \de (2\theta+n(1-\theta))}{\theta} \right]+1,
		\end{equation}
		which does not depend on $\eps$. Using  equation \eqref{eq: tau4}, it is an elementary exercise to verify that \eqref{eq: M_n size} holds for the choice of $M$ in \eqref{eq: M_n choice}, and the proof is complete.
	\end{proof}
	
	For each $k\in[n_0, n]$ and any color $\xi\in \Xi$, fix the mappings $F^\xi_k$, as in \cref{lem:
		choice of v_j}. Moreover, define the mappings $\widetilde{F}^\xi_k$ similarly to $F^\xi_k$, but using
	the corresponding decreasing order on $J_k(\xi)$ to pick the vectors $\widetilde{v}_j=\widetilde{v}_j^k$
	inductively as in the proof of \cref{lem: choice of v_j}. Define the map $F_k : X \to
	\mathbb{R}^{\widetilde{N}}$, where $\widetilde{N} = 2N_n M_n$ and $N_n=|\Xi|$, to be the direct
	product of all maps $F^\xi_k$ and $\widetilde{F}^\xi_k$ for all colors $\xi \in \Xi$. We prove that the $n$-th such map, i.e., $$
	F:=F_n=\bigotimes_{\xi\in\Xi}(F^\xi\otimes \widetilde{F}^\xi): X\rightarrow \R^{\widetilde{N}},
	$$
	%
	satisfies the conditions of the weak embedding in the statement of \cref{thm:main}.
	
	\begin{proof}[Proof of \cref{thm:main}]
		Given $x, y \in X$, fix $k\in [n_0,n]$ to be the largest integer such that
		\begin{equation}\label{eq: d(x,y) smaller than 4r_k}
			d(x, y) \leq 4r_{k-1} = 4\tau^{-2}r_k.
		\end{equation} Fix $\tau$ and $F$ as above. By the Lipschitz norm estimate of $F_k^\xi$ in \eqref{eq: F_k^ksi Lip norm} we have
		$$
		|F^\xi_k(x) - F^\xi_k(y)| \leq ||F^\xi_k||_{\text{Lip}}\,d(x, y) \leq r_k^{\eps-1}d(x, y),
		$$ which by \eqref{eq: d(x,y) smaller than 4r_k} implies
		\begin{equation}\label{eq: F_k upper Lip}
			|F^\xi_k(x) - F^\xi_k(y)|  \leq \left(d(x, y) / 4\tau^{-2}\right)^{\eps-1}d(x, y) = (4\tau^{-2})^{1-\eps}d(x, y)^\eps.
		\end{equation}
		
		Due to \eqref{eq: F_k^ksi convergence} we also have
		\begin{equation}\label{eq: F_k lower Lip}
			\left| |F^\xi(x) - F^\xi(y)| - |F^\xi_k(x) - F^\xi_k(y)| \right| \leq 2||F^\xi - F^\xi_k||_\infty \leq 4\tau^2r^\eps_{k+1} = 4 \tau^2\tau^{2\eps}r_k^\eps\leq 2\tau^2 d(x,y)^\eps
		\end{equation}
		Hence, by \eqref{eq: F_k upper Lip}, \eqref{eq: F_k lower Lip}, and similar estimates for $\widetilde{F}_\xi^k$, the Lipschitz estimate 
		$$
		|F(x) - F(y)| \leq 5\sqrt{N_n} \tau^{-2(1-\eps)}d(x, y)^\eps
		$$
		follows by definition of $F$ as the direct product of these $N_n$ maps.
		
		For the lower bound further assume that
		\begin{equation}\label{eq: d(x,y) not too small condition}
			4r_k \leq d(x, y) \leq 4r_{k-1} = 4\tau^{-2}r_k.
		\end{equation}
		By \eqref{eq:X_covering_def} there is $j \in J_k$ such that $x \in B_j$. Fix the color $\xi \in \Xi$ such that $j \in J_k(\xi)$.

		By $d(x, x_j) \leq r_k$, since $x \in B_j$, and \eqref{eq: d(x,y) not too small condition}, we have that $y \in B(x_j, 10\tau^{-2}r_k)$. On the other hand, if $y \in 2B_j$, then $d(x, y) \leq d(x, x_j) + d(x_j, y) \leq 3r_k$, which contradicts \eqref{eq: d(x,y) not too small condition}. Hence, we have shown that
		\begin{equation}\label{eq: y is out of 2B}
			y \in B(x_j, 10\tau^{-2}r_k) \setminus 2B_j,
		\end{equation} which allows to apply \cref{lem: choice of v_j} in what follows.
		
		Suppose $y \in 2B_i$ for some $i \in J_k(\xi)$, with $i \neq j$ by \eqref{eq: y is out of 2B}.
		Assume that $i < j$, and the proof in the case $i>j$ follows similarly by considering the map $\widetilde{F}^\xi_k$ instead of $F^\xi_k$ in what follows. Recall that  $\phi_l(y)=0$ for all $l \neq i$, by \eqref{eq:phi_0_outside_2B}, \eqref{eq:same_color_points_far}, and \eqref{eq: tau1}. Then by \eqref{eq: Fk is F_k-1 + phis}, \eqref{eq: G_k def} we have
		\begin{equation}\label{eq: F_k and G_k}
			F^\xi_k(y) = F^\xi_{k-1}(y) + r^\eps_k v_i \phi_i(y) = G^\xi_{k,j}(y),
		\end{equation}
		with the convention that $F^\xi_{k-1}(y) = 0$ for $k = n_0$. By \eqref{eq: y is out of 2B}, we can
		apply \cref{lem: choice of v_j}, and \eqref{eq: Lower Lip Lemma} implies that
		$$
		|F^\xi_k(x) - F^\xi_k(y)| = |F^\xi_k(x) - G^\xi_{k,j}(y)| \geq \tau^3 r^\eps_k.
		$$
		By the above relation and the second to last inequality in \eqref{eq: F_k lower Lip} and get that
		$$
		|F^\xi(x) - F^\xi(y)| \geq |F^\xi_k(x) - F^\xi_k(y)| - 4\tau^2 \tau^{2\eps} r^\eps_k\geq \tau^3 r_k^\eps-4\tau^2r_{k+1}^\eps.
		$$
		However, we have $\tau^3 r_k^\eps(1-4\tau^{-1}\tau^{2\eps})\geq \tau^3 r^\eps_k/2$ by \eqref{eq:rk_scales_def} and \eqref{eq: tau3}, which gives
		$$
		|F^\xi(x) - F^\xi(y)| \geq \frac{\tau^3 r_k^\eps}{2}.
		$$
		Using the above on the trivial inequality $|F(x) - F(y)| \geq |F^\xi(x) - F^\xi(y)|$, along with \eqref{eq: d(x,y) not too small condition}, implies
		\begin{equation}\label{eq: F lower-Lip}
			|F(x) - F(y)| \geq \frac{\tau^3}{2}\left( \frac{d(x,y)}{4\tau^{-2}} \right)^\eps \geq \frac{\tau^{5}}{8} d(x, y)^\eps.
		\end{equation}
		The proof of the lower estimate is thus complete for the case where $y \in 2B_i$ for some $i \in J_k(\xi)$. 
		
		Suppose that $y \notin 2B_i$ for all $i \in J_k(\xi)$. Then $\phi_i(y)=0$ for all $i\in J_k(\xi)$ by \eqref{eq:phi_0_outside_2B}. This implies that $F^\xi_k(y) = F^\xi_{k-1}(y)= G^\xi_{k,j}(y)$ for $j$ with $x \in B_j$. Thus, \eqref{eq: F_k and G_k} holds and the proof follows identically to the previous case. As a result, the map $F$ satisfies all the necessary properties and the proof is complete.
	\end{proof}

	\section{Final remarks}\label{sec:questions}
	The proof of \cref{thm:main} relied on properties of the Assouad spectrum, namely the quasidoubling
	property of the space. It would be of independent interest to construct such an embedding using
	solely the fact that the Minkowski dimension is finite, without employing the Assouad spectrum.
	Methods that depend on colorings of nets to define the embedding, such as those in the proof of
	Assouad's theorem or \cref{thm:main}, would no longer be applicable directly.
	
	In addition, Assouad's theorem is in fact a sharp result, in the sense that if $(X,d)$ is a metric space that can be snowflaked and bi-Lipschitz embedded in a Euclidean space, then it has finite Assouad dimension (see \cite{hei:lectures}). It is not clear whether such a ``converse" also applies in the case of the Minkowski dimension. Namely, if there is a sequence of such weak embeddings $F_n$ of $X$ into a sequence of Euclidean spaces $\R^{M_n}$ of increasing dimensions $M_n$ as in \cref{thm:main}, it is not clear whether the Minkowski dimension of $X$ needs to be finite.
	
	Moreover, it seems that embedding properties of metric spaces are in close relation with certain
	dimension notions. As mentioned in the introduction, G.C. David \cite{Guy_Nagata} proved a similar
	weak embedding result for spaces with finite Nagata dimension. That dimension notion, while quite
	close to an integer version of the Assouad dimension, has more of a topological character. It is
	reasonable to ask whether other fractal (and fractional) dimensions could lead to such embedding
	results for metric spaces, such as the finiteness of the Hausdorff dimension of a space. We
	emphasize that a direct proof of \cref{thm:main} from the finitenes of $\ovdimB X$, without using the finiteness of the Assouad spectrum, could possibly shed light on the Hausdorff dimension case as well. That might be possible by employing the Intermediate Dimensions and their properties, as introduced by Falconer, Fraser and Kempton \cite{Interm_dim_Falc_Fras_kem}.

	\bibliographystyle{acm}

\end{document}